\documentclass{psp}
\usepackage[latin1]{inputenc}
\usepackage[T1]{fontenc}
\usepackage{parskip}
\usepackage{amsmath,amsfonts,amssymb,amsxtra}
\usepackage{latexsym}
\usepackage{theorem}
\usepackage{graphicx,psfrag}
\usepackage{verbatim}

{\theoremstyle{plain}                       
\theorembodyfont{\itshape}
\newtheorem{Theorem}{Theorem}}

{\theoremstyle{plain}                       
\theorembodyfont{\itshape}
\newtheorem{Corollary}{Corollary}}

{\theoremstyle{plain}                       
\theorembodyfont{\itshape}
\newtheorem{Proposition}{Proposition}}

{\theoremstyle{plain}                       
\theorembodyfont{\rmfamily}
\newtheorem{Definition}{Definition}}

{\theoremstyle{plain}                       
\theorembodyfont{\rmfamily}
\newtheorem{Notation}{Notation}}

{\theoremstyle{plain}                       
\theorembodyfont{\rmfamily}
\newtheorem{Example}{Example}}

{\theoremstyle{plain}                       
\theorembodyfont{\rmfamily}
\newtheorem{Remark}{Remark}}

{\theoremstyle{plain}                       
\theorembodyfont{\itshape}
\newtheorem{Lemma}{Lemma}}

\begin{document}

\title[Coding map for a contractive Markov system]{Coding map for a contractive Markov system}
\author[Ivan Werner]{Ivan Werner\thanks{Supported by EPSRC and School of
Mathematics and Statistics of University of St Andrews.}\\
{\small\it Mathematical Institute,
   University of St Andrews,
   St Andrews, Fife KY16 9SS, UK}\\
   {\small\it e-mail: ivan\_werner@pochta.ru}}
\maketitle

\begin{abstract}\noindent
 In this paper, we develop the theory of
 contractive Markov systems initiated in \cite{Wer1}.
 We construct a coding map for such systems and investigate some of its properties. \\

\end{abstract}

\section{Introduction}

\noindent In  \cite{Wer1}, we introduced a theory of {\it
contractive Markov systems (CMSs)} which provides a unifying
framework in 'fractal' geometry. It extends the known theory of {\it
iterated function systems (IFSs) with place dependent
probabilities}, which are contractive on average,
\cite{BDEG}\cite{Elton} in a way that it also covers {\it graph
directed constructions} of 'fractal' sets \cite{MW}. This theory
also is a natural generalization of the theory of finite Markov
chains.

\noindent A contractive Markov system  with an average contracting
rate $0<a<1$ is a family
\[\left(K_{i(e)},w_e,p_e\right)_{e\in E}\]
(see Fig. 1) where $E$ is the set of edges of a directed
(multi)graph $(V,E,i,t)$ ($V:=\{1,...,N\}$ is the set of vertices
of the directed (multi)graph (we do not exclude the case $N=1$),
$i:E\longrightarrow V$ is a map indicating the initial vertex of
each edge and $t:E\longrightarrow V$ is a map indicating the
terminal vertex of each edge), $K_1,K_2,...,K_N$ is a partition of
a metric space $(K,d)$ into non-empty Borel subsets, $(w_e)_{e\in
E}$ is a family of Borel measurable self-maps on the metric space
such that $w_e\left(K_{i(e)}\right)\subset K_{t(e)}$ for all $e\in
E$ and $(p_e)_{e\in E}$ is a family of Borel measurable
probability functions on $K$ (i.e. $p_e(x)\geq 0$ for all $e\in E$
and $\sum_{e\in E}p_e(x)=1$ for all $x\in K$) (associated with the
maps) such that each $p_e$ is zero on the complement of $K_{i(e)}$
and the system satisfies the following {\it condition of a
contractiveness on average}
\begin{equation}\label{cc}
 \sum\limits_{e\in E}p_e(x)d(w_ex,w_ey)\leq ad(x,y)\mbox{ for all
}x,y\in K_i,\ i=1,...,N.
\end{equation}

This condition was discovered by R. Isaac \cite{Is} for the case
$N=1$.

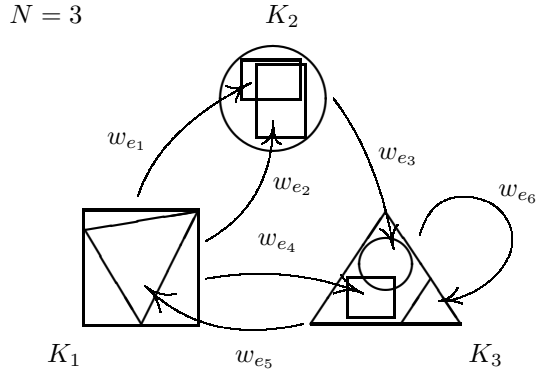
\begin{center}
\unitlength 1mm
\begin{picture}(70,70)\thicklines
\put(35,50){\circle{20}} \put(10,20){\framebox(15,15)}
\put(40,20){\line(2,3){10}} \put(40,20){\line(4,0){20}}
\put(50,35){\line(2,-3){10}} \put(5,15){$K_1$} \put(34,60){$K_2$}
\put(61,15){$K_3$} \put(31,50){\framebox(7.5,5)}
\put(33,45){\framebox(6.25,9.37)} \put(50,28){\circle{7.5}}
\put(45,21){\framebox(6,5)} \put(10,32.5){\line(6,1){15}}
\put(10,32.5){\line(3,-5){7.5}} \put(17.5,20){\line(1,2){7.5}}
\put(52,20){\line(2,3){4}} \put(13,44){$w_{e_1}$}
\put(35,38){$w_{e_2}$} \put(49,42){$w_{e_3}$}
\put(33,30.5){$w_{e_4}$} \put(30,15){$w_{e_5}$}
\put(65,37){$w_{e_6}$} \put(15,5){ Fig. 1. A Markov system.}
\put(0,60){$N=3$} \thinlines \linethickness{0.1mm}
\bezier{300}(17,37)(20,46)(32,52)
\bezier{50}(32,52)(30.5,51.7)(30,49.5)
\bezier{50}(32,52)(30,51)(28.7,51.7)
\bezier{300}(26,31)(35,36)(35,47)
\bezier{50}(35,47)(35,44.5)(33.5,44)
\bezier{50}(35,47)(35,44)(36,44) \bezier{300}(43,50)(49,42)(51,30)
\bezier{50}(51,30)(50.5,32)(49.2,32.6)
\bezier{50}(51,30)(50.6,32)(51.5,33.2)
\bezier{300}(39,20)(26,17)(18,25)
\bezier{50}(18,25)(19.5,24)(20,21.55)
\bezier{50}(18,25)(20,23.5)(22,24)
\bezier{300}(26,26)(37,28)(47,24) \bezier{50}(47,24)(45,25)(43,24)
\bezier{50}(47,24)(45,25)(44,26.5)
\bezier{100}(54.5,31.9)(56,37.3)(61,36.9)
\bezier{100}(61,36.9)(64.5,36.5)(66,34)
\bezier{100}(66,34)(68,30.5)(64.9,26.8)
\bezier{100}(64.9,26.8)(61.6,23.3)(57,23)
\bezier{50}(57,23)(58.5,23.3)(60.1,22.7)
\bezier{50}(57,23)(58.8,23.3)(59.5,24.8)
\end{picture}
\end{center}

Such a system determines a Markov operator $U$ on the set of all
bounded Borel measurable functions $\mathcal{L}^0(K)$ by
\[Uf:=\sum\limits_{e\in E}p_ef\circ w_e\mbox{ for all
 }f\in\mathcal{L}^0(K)\] and its adjoint operator $U^*$ on the set of all Borel probability
 measures $P(K)$ by
\[U^*\nu(f):=\int U(f)d\nu=\sum\limits_{e\in E}\int
\limits_{K_{i(e)}}p_ef\circ w_ed\nu\mbox{ for all
}f\in\mathcal{L}^0(K)\mbox{ and }\nu\in P(K).\]

\begin{Remark}
    Note that each maps $w_e$ and each probability function $p_e$ need to
    be defined only on the corresponding vertex set $K_{i(e)}$. This
   is sufficient for the condition (\ref{cc}) and the definition of
   $U^*$. For the definition of $U$, we can consider each $w_e|_{K_{i(e)}}$ to
   be extended on the whole space $K$ arbitrarily and each $p_e|_{K_{i(e)}}$ to
   be extended on $K$ by zero.

   Also, the situation applies where each vertex set $K_i$ has its
   own metric $d_i$. In this case, one can set
   \[d(x,y)=\left\{\begin{array}{cc}
       d_i(x,y) & \mbox{ if }x,y\in K_i \\
       \infty & \mbox{otherwise}
     \end{array}\right.
   \] and use the convention $0\times\infty=0$.
\end{Remark}

We say $\mu\in P(K)$ is an {\it invariant measure} of the CMS iff
$U^*\mu=\mu$. A Borel probability measure $\mu$ is called {\it
 attractive} measure of the CMS if
 \[{U^*}^n\nu\stackrel{w^*}{\to}\mu\mbox{ for all }\nu\in P(K),\]
 where $w^*$ means weak$^*$ convergence.
Note that an attractive Borel probability measure is a unique
invariant Borel probability measure of the CMS if $U$ maps
continuous functions on continuous functions.

It was shown in \cite{Wer1} that operator $U^*$ inherits some
properties from its trivial case of a finite Markov chain if the
vertex sets $K_1,...,K_N$ form an open partition of the state
space and the restrictions of the probability functions on their
vertex sets are Dini-continuous and bounded away from zero.
Namely, it has a unique invariant probability measure in an
irreducible case and an attractive probability measure in an
aperiodic case.

The coding map is an important tool in 'fractal' geometry which
allows one to represent a constructed set as an image of a code
space under this map, that is, to code elements in this set by
infinite sequences of elements of $E$. We shall denote
$\Sigma:=E^{\mathbb{Z}}$. Let's see some examples.

\begin{Example}[decimal expansion]
 Consider ten maps $w_e$, $e\in E:=\{0,...,9\}$, on $([0,1], |.|)$ given by
 $w_e(x):=1/10x+e/10$ for all $x\in [0,1]$. Obviously, for any
 family of probability functions $p_e$, $e\in E$, the family
 $([0,1],w_e,p_e)_{e\in E}$ is a CMS. Fix an $x\in
 [0,1]$. We define the coding map
 $F:\Sigma \longrightarrow\ [0,1]$ by
 \begin{equation}\label{cm0}
 F(\sigma):=\lim\limits_{m\to-\infty}w_{\sigma_0}\circ...\circ
w_{\sigma_m}(x)\mbox{ for
 all }\sigma\in\Sigma,
\end{equation}
 which is nothing else as the decimal expansion of real numbers from
 $[0,1]$.
\end{Example}

\begin{Example}[Cantor set]
 Consider two maps  $w_e$, $e\in E:=\{0,1\}$, on $([0,1], |.|)$ given by
 $w_e(x):=1/3x+e2/3$ for all $x\in [0,1]$. Again, for any
 family of probability functions $p_e$, $e\in E$, the family
 $([0,1],w_e,p_e)_{e\in E}$ is a CMS. Fix an $x\in
 [0,1]$. Then the coding map
 $F:\Sigma \longrightarrow\ [0,1]$ given by the limit (\ref{cm0}), for
 all $\sigma\in\Sigma$, is nothing else as the binary expansion of the Cantor set.
\end{Example}

\begin{Example}
 Consider two maps  $w_e$, $e\in E:=\{0,1\}$, on $(\mathbb{R}, |.|)$ given by
 $w_0(x):=1/2x$ and $w_1(x):=-2x+3$ for all $x\in
 \mathbb{R}$ with constant probability functions $p_0:=3/4$ and $p_1:=1/4$. A simple calculation
 shows that
 $([0,1],w_e,p_e)_{e\in E}$ defines a CMS with an average contraction rate $7/8$. Fix an $x\in
 \mathbb{R}$. It was shown in \cite{BE} that the coding map
 $F:\Sigma \longrightarrow \mathbb{R}$ given by the limit (\ref{cm0})
  exists for $P$-a.e. $\sigma\in\Sigma$ and is independent of the choice of $x$ modulo a $P$-zero set
  where $P$ is a Bernoulli measure on $\Sigma$ given by
  $P([e_1,e_2...,e_k])=p_{e_1}p_{e_2}...p_{e_k}$ for every thin cylinder set
  $[e_1,...,e_k]\subset\Sigma$.
\end{Example}

For a CMS $\left(K_{i(e)},w_e,p_e\right)_{e\in E}$ with several
vertex sets (as at Fig. 1), we are going to define the coding map
as follows. For each $i$, fix $x_i\in K_i$ and set
\begin{equation}\label{cm}
    F(\sigma):=\lim\limits_{m\to-\infty}w_{\sigma_0}\circ...\circ
w_{\sigma_m}(x_{i(\sigma_m)})
\end{equation} for $\sigma\in\Sigma$ if the limit
exists. It is a simple exercise for the reader to show that
$F(\sigma)$ exists for all $\sigma\in\Sigma_G:=\{\sigma\in\Sigma:\
t(\sigma_m)=i(\sigma_{m-1})\ \forall m\in\mathbb{Z}\}$, it is
independent of the choice of $x_i$'s and $F$ is
H\"{o}lder-continuous if all maps $w_e|_{K_{i(e)}}$ are contractions
(see \cite{K} for the case $N=1$).

We would like to illustrate it by the following example, which is
associated with the notion of $g$-measures (see e.g. \cite{Wal})
\begin{Example}
  Let $G:=(V,E,i,t)$ be a finite irreducible directed
  (multi)graph. Let $\Sigma^-_G:=\{(...,\sigma_{-1},\sigma_0):\
  \sigma_m\in E\mbox{ and }
t(\sigma_m)=i(\sigma_{m-1})\ \forall
m\in\mathbb{Z}\setminus\mathbb{N}\}$ ({\it be one-sided
  subshift of finite type} associated with $G$) endowed with the
  metric $d(\sigma,\sigma'):=2^k$ where $k$ is the smallest
  integer with $\sigma_i=\sigma'_i$ for all $k<i\leq 0$. Let $g$
  be a positive, Dini-continuous function on $\Sigma_G$ such that
  \[\sum\limits_{y\in T^{-1}(\{x\})}g(y)=1\mbox{ for all }x\in\Sigma_G\]
  where $T$ is the right shift map on $\Sigma^-_G$. Set
  $K_i:=\left\{\sigma\in\Sigma^-_G:t(\sigma_0)=i\right\}$  for
  every $i\in V$ and, for
  every $e\in E$,
  \[w_e(\sigma):=(...,\sigma_{-1},\sigma_{0},e),\ p_e(\sigma):=g(...,\sigma_{-1},\sigma_{0},e)
  \mbox{ for all }\sigma\in K_{i(e)}.\]
  Obviously, maps $(w_e)_{e\in E}$ are contractions. Therefore,
  $\left(K_{i(e)}, w_e, p_e\right)_{e\in E}$
   defines a CMS. In this example, the coding
  map $F:\Sigma_G\longrightarrow\ \Sigma^-_G$ given by (\ref{cm})
  is nothing else as the natural projection.
\end{Example}

In this paper, we are concerned with the question whether limit
(\ref{cm}) exists almost everywhere  and is independent of the
choice of $x_i$ up to a set of measure zero with respect to a
natural measure if the maps $w_e$ are contractive only on average
and the probabilities $p_e|_{K_{i(e)}}$ are place-dependent (see the
next example). We show that this is true under some conditions on
the probability functions.

\begin{Example}
 Let $\mathbb{R}^2$ be normed by $\|.\|_1$. Let
 $K_1:=\{(x,y)\in\mathbb{R}^2:\ y\geq 1/2\}$ and $K_2:=\{(x,y)\in\mathbb{R}^2:\ y\leq
 -1/2\}$. Consider the following maps on $\mathbb{R}^2$:
 \begin{eqnarray*}
   &&w_1\left(\begin{array}{c}x\\ y\end{array}\right):=\left(\begin{array}{c}
   -\frac{1}{2}x-1\\-\frac{3}{2}y+\frac{1}{4}\end{array}\right),\
   w_2\left(\begin{array}{c}x\\ y\end{array}\right):=\left(\begin{array}{c}
   -\frac{3}{2}x+1\\\frac{1}{4}y+\frac{3}{8}\end{array}\right),\mbox{
   and }\\
  &&w_3\left(\begin{array}{c}x\\ y\end{array}\right):=\left(\begin{array}{c}
   -\frac{1}{2}|x|+1\\ -\frac{3}{4}y+\frac{1}{8}\end{array}\right)
 \end{eqnarray*}
 with probability functions
 \begin{eqnarray*}
  && p_1\left(\begin{array}{c}x\\ y\end{array}\right):=\left(\frac{1}{15}\sin^2\|(x,y)\|_1+\frac{53}{105}
   \right)1_{K_1}(x,y),\\
  && p_2\left(\begin{array}{c}x\\ y\end{array}\right)
   :=\left(\frac{1}{15}\cos^2\|(x,y)\|_1+\frac{3}{7}\right) 1_{K_1}(x,y)\mbox{ and }
   p_3:=1_{K_2}.
 \end{eqnarray*}
 A simple calculation shows that $(K_{i(e)},w_e,p_e)_{e\in\{1,2,3\}}$, where $i(1)=1$, $i(2)=1$
 and $i(3)=2$, defines a CMS with an average
 contracting rate $209/210$. (By Theorem 2 in \cite{Wer1}, it
 has a unique (attractive) invariant probability measure.)
\end{Example}

Our proof consists of two parts. In the first part, we construct a
suitable outer measure on the code space and show that, for every
CMS, limit (\ref{cm}) exists almost everywhere with respect to
this outer measure and is independent of the choice of $x_i$'s up
to a set of measure zero. It involves some measure-theoretic
technique which seems to be new. In the second part, we show that
the natural measure on the code space, which we call the {\it
generalized Markov measure}, is absolutely continuous with respect
to the outer measure if the restrictions of probability functions
on their vertex sets are Dini-continuous and bounded away from
zero. This implies the desired result.

The coding map which we construct here opens a way for various
applications. In \cite{Wer4}, we compute the Kolmogorov-Sinai
entropy $h_M(S)$ of the generalized Markov shift associated with the
CMS (see Definition \ref{gMs}) using the coding map. We show that
\[h_M(S)=\sum\limits_{e\in E}\int
\limits_{K_{i(e)}}p_e\log p_ed\mu,\] where $\mu$ is a unique
invariant Borel probability measure of the CMS. In \cite{Wer5}, we
prove an ergodic theorem for CMSs using the coding map. In
\cite{Wer6}, we show that the generalized Markov measure is a unique
equilibrium state with respect to an energy function the
construction of which involves the coding map.

\section{Construction with respect to an outer measure}

Let $\left(K_{i(e)},w_e,p_e\right)_{e\in E}$ be a contractive
Markov system with an average contracting rate $0<a<1$ and an
invariant Borel probability measure $\mu$. We assume metric space
$(K,d)$ to be complete, the set of edges $E$ to be finite. We do
not pose any conditions on the directed graph.

Let $\Sigma:=\{(...,\sigma_{-1},\sigma_0,\sigma_1,...):e_i\in E\
\forall i\in\mathbb{Z}\}$. We shall denote by $d'$ the metric on
$\Sigma$ defined by $d'(\sigma,\sigma'):=(1/2)^k$ where $k$ is the
largest integer with $\sigma_i=\sigma'_i$ for all $|i|<k$. Let $S$
be the left shift map on $\Sigma$. We call the set $
_m[e_m,...,e_n]:=\{\sigma\in\Sigma:\ \sigma_i=e_i\mbox{ for all
}m\leq i\leq n\}$,for $m\leq n\in\mathbb{Z}$, a {\it cylinder}.
Denote by $\mathcal{A}$ the finite $\sigma$-algebra generated by the
zero time partition $\{_0[e]:e\in E\}$ of $\Sigma$, and define, for
each integer $m\leq 1$,
\[\mathcal{A}_m:=\bigvee\limits_{i=m}^{+\infty} S^{-i}\mathcal{A},\]
which is the smallest $\sigma$-algebra containing  all finite
 $\sigma$-algebras $\bigvee_{i=m}^{n}
S^{-i}\mathcal{A}$, $n\geq m$. Let $x\in K$. For each integer
$m\leq 1$ let $P_x^m$ be the probability measure on the
$\sigma$-algebra $\mathcal{A}_m$ given by
\[P^m_x( _{m}[e_{m},...,e_n])=p_{e_{m}}(x)p_{e_{m+1}}(w_{e_{m}}(x))...p_{e_n}(w_{e_{n-1}}\circ...\circ
w_{e_{m}}(x))\] for all cylinder sets $_{m}[e_{m},...,e_n]$,
$n\geq{m}$.
\begin{Lemma}
 Let $m\leq 1$ and $A\in\mathcal{A}_m$.
 Then $x\longmapsto P_x^m(A)$ is a Borel measurable function on $K$.
\end{Lemma}
\begin{proof} Set
\[\mathcal{D}:=\left\{A\in\mathcal{A}_m:K\ni x\longmapsto P_x^m(A)\mbox{
Borel measurable }\right\}.\] Then, by definition of $P_x^m$,
$\mathcal{D}$ contains all cylinders of the form
$_{m}[e_{m},...,e_n],\ n\geq m$, which generate $\mathcal{A}_m$.
Furthermore, obviously it holds true that
\[\Sigma\in\mathcal{D},\]
\[A\in\mathcal{D}\Rightarrow\Sigma\setminus A\in\mathcal{D}\]
and, for any pairwise disjoint family
$(A_n)_{n\in\mathbb{N}}\subset\mathcal{D}$,
\[\bigcup\limits_{n\in\mathbb{N}}A_n\in\mathcal{D},\]
i.e. $\mathcal{D}$ is a Dynkin-system. Hence, $\mathcal{D}$ contains
the Dynkin-system which is generated by the cylinders. Since the set
of the cylinders is $\cap$-stable, it follows that
$\mathcal{D}=\mathcal{A}_m$.\end{proof}

\begin{Definition}
 Let $\nu\in P(K)$. We call a probability measure $\Phi_m(\nu)$ on
 $(\Sigma,\mathcal{A}_m)$ given by
 \[\Phi_m(\nu)(A):=\int P_x^m(A)d\nu(x),\
 A\in\mathcal{A}_m,\] {\it the $m$-th lift of} $\nu$.
\end{Definition}
\begin{Definition}
 Set
 \[\mathcal{C}(B):=\left\{(A_m)_{m=0}^{-\infty}:A_m\in\mathcal{A}_m\
 \forall m\mbox{ and }B\subset\bigcup\limits_{m=0}^{-\infty}A_m
 \right\}\] for $B\subset\Sigma$. Let $\nu\in P(K)$. We call a set
 function given by
 \[\Phi(\nu)(B):=\inf\left\{\sum\limits_{m=0}^{-\infty}\Phi_m(\nu)(A_m):(A_m)_{m\leq 0}\in\mathcal{C}(B)\right\}, B\subset\Sigma,\]
{ \it  the lift of} $\nu$.
\end{Definition}
\begin{Lemma}\label{om}
 Let $\nu,\lambda\in P(K)$. Then\\
 $(i)$ $\Phi(\nu)$ is an outer measure on $\Sigma$.\\
 $(ii)$ If $\Phi_m(\nu)\ll\Phi_m(\lambda)$ for all $m\leq 0$ , then for all $\epsilon>0$
 there exists $\delta>0$ such that
 \[\Phi(\lambda)(B)<\delta\ \Rightarrow\ \Phi(\nu)(B)<\epsilon\mbox{ for all
 }B\subset\Sigma.\]
\end{Lemma}
\begin{proof} It is obvious that $\Phi(\nu)(\emptyset)=0$.

 Let
$B_1\subset B_2\subset\Sigma$. Then
$\mathcal{C}(B_1)\supset\mathcal{C}(B_2)$ and therefore
\[\Phi(\nu)(B_1)\leq\Phi(\nu)(B_2).\]

Now, we show
\[\Phi(\nu)\left(\bigcup\limits_{i=1}^{\infty}B_i\right)\leq\sum\limits_{i=1}^\infty\Phi(\nu)(B_i)\]
for all $B_i\subset\Sigma$, $i\in\mathbb{N}$. We can assume that
the right hand side is finite. Let $\epsilon>0$. Then for every
$i\in\mathbb{N}$ there exists $(A_{im})_{m\leq
0}\in\mathcal{C}(B_i)$ such that
\[\Phi(\nu)(B_i)>\sum\limits_{m=0}^{-\infty}\Phi_m(\nu)(A_{im})-\epsilon
2^{-i}.\] Since $\left(\bigcup_{i=1}^\infty A_{im}\right)_{m\leq
0}\in\mathcal{C}\left(\bigcup_{i=1}^{\infty}B_i\right)$, it
follows that
\begin{eqnarray*}
 \Phi(\nu)\left(\bigcup\limits_{i=1}^{\infty}B_i\right)&\leq&\sum\limits_{m=0}^{-\infty}\Phi_m(\nu)\left(\bigcup_{i=1}^\infty
A_{im}\right)\\
&\leq&\sum\limits_{m=0}^{-\infty}\sum\limits_{i=1}^\infty\Phi_m(\nu)\left(
A_{im}\right)\\
&\leq&\sum\limits_{i=1}^\infty\Phi(\nu)(B_i)+\epsilon,
\end{eqnarray*}
as desired.

 Suppose the claim in (ii) is not true.
Then there exist $\epsilon>0$ and a sequence of subsets
$(B_n)_{n\in\mathbb{N}}$ such that
\[\Phi(\lambda)(B_n)<2^{-n}\mbox{ and
}\Phi(\nu)(B_n)\geq\epsilon\mbox{ for all }n\in\mathbb{N}.\]
 Then
for every $n\in\mathbb{N}$ there exist $(A_{nm})_{m\leq
0}\in\mathcal{C}(B_n)$ such that
\[\sum\limits_{m=0}^{-\infty}\Phi_m(\lambda)(A_{nm})<2^{-n}.\]
Set \[D_m:=\bigcap\limits_{k=1}^\infty\bigcup\limits_{n\geq
k}A_{nm}\mbox{ for each }m\leq 0.\] Then, for each $m\leq 0$,
\[\Phi_m(\lambda)(D_m)\leq \Phi_m(\lambda)\left(\bigcup\limits_{n\geq
k}A_{nm}\right)\leq \sum\limits_{n\geq k}2^{-n}\] for all $k\geq
0$. Hence $\Phi_m(\lambda)(D_m)=0$ for all $m\leq 0$. Since
$\Phi_m(\nu)\ll\Phi_m(\lambda)$ for all $m\leq 0$, this implies
that
\[0=\Phi_m(\nu)(D_m)=\lim\limits_{k\to\infty}\Phi_m(\nu)\left(\bigcup\limits_{n\geq
k}A_{nm}\right)\geq
\limsup\limits_{k\to\infty}\Phi_m(\nu)\left(A_{km}\right)\] for
all $m\leq 0$. Hence
\[\limsup\limits_{k\to\infty}\Phi(\nu)(B_k)\leq\limsup\limits_{k\to\infty}\sum\limits_{m=0}^{-\infty}\Phi_m(\nu)\left(A_{km}\right)=0,\]
which is a contradiction.
 \end{proof}

  We use further the following notation.
 \begin{Notation} Fix $x_i\in K_i$ for each $i\in\{1,...,N\}$ and
 set
 \[P^m_{x_1...x_N}:=\Phi_m\left(\frac{1}{N}\sum\limits_{i=1}^N\delta_{x_i}\right)\mbox{ and }
 P_{x_1...x_N}:=\Phi\left(\frac{1}{N}\sum\limits_{i=1}^N\delta_{x_i}\right)\]
 \end{Notation}
for every $m\in\mathbb{Z}\setminus\mathbb{N}$, where $\delta_x$
denotes the Dirac probability measure concentrated at $x$. Then
\[P^m_{x_1...x_N}( _m[e_m,...,e_n])=\frac{1}{N}p_{e_m}(x_{i(e_m)})p_{e_{m+1}}(w_{e_m}x_{i(e_m)})...p_{e_n}(w_{e_{n-1}}
\circ...\circ w_{e_m}x_{i(e_m)})\] for every cylinder set
$_m[e_m,...,e_n]$.

 Now, for every  $m\leq 0$ and $n\geq m$ define a  random variable
  \begin{eqnarray*}
    Y_{mn}^{x_1...x_N}:\Sigma&\longrightarrow& K  \\
                     \sigma&\longmapsto& w_{\sigma_n}\circ
                     w_{\sigma_{n-1}}\circ...\circ w_{\sigma_{m}}(x_{i(\sigma_{m})})
  \end{eqnarray*}
with respect to the measure $P^m_{x_1...x_N}$.

Now, we are going to prove the main lemma which enables us to define
the coding map. The proof of it involves a more general version of
Borel-Cantelli argument than that which was used by Barnsley and
Elton in \cite{BE} (they considered the case $N=1$ with constant
probabilities). Their key point, the reversion of the order of
finite sequences of the maps, does not work here because
$P^m_{x_1...x_N}( _m[e_m,e_{m+1}...,e_n])\neq P^m_{x_1...x_N}(
_m[e_n,e_{n-1}...,e_m])$ in general. That is why we first needed to
construct the outer measure $P_{x_1...x_N}$. The kind of
Borel-Cantelli argument which we apply here seems to be new in its
generality.

\begin{Lemma}\label{cml}
  Let $x_i,y_i \in K_i$ for each $1\leq i\leq N$.\\
(i)
\[\lim\limits_{m\to-\infty}d\left(Y^{x_1...x_N}_{m0},Y^{y_1...y_N}_{m0}\right)=0\
P_{x_1...x_N}\mbox{-a.e.},\]
(ii)
\[F_{x_1...x_N}:=\lim_{m\to-\infty}Y^{x_1...x_N}_{m0}\mbox{ exists }
P_{x_1...x_N}\mbox{-a.e.},\] and by $(i)$
$F_{x_1...x_N}=F_{y_1...y_N}$ $P_{x_1...x_N}$-a.e..

 (iii) There exists a sequence of closed subsets $Q_1\subset
 Q_2\subset...\subset\Sigma$ with \\ $\sum\limits_{k=1}^\infty P_{x_1...x_N}(\Sigma\setminus Q_k)<\infty$
such that $F_{x_1...x_N}|_{Q_k}$ is locally H\"{o}lder-continuous
with the same H\"{o}lder-constants for all $k\in\mathbb{N}$, i.e.
there exist $\alpha, C>0$ such that for every $k$ there exists
$\delta_k>0$ such that
\[ \sigma,\sigma'\in Q_k\mbox{ with }d'(\sigma,\sigma')\leq\delta_k\
\Rightarrow\ d(F_{x_1...x_N}(\sigma),F_{x_1...x_N}(\sigma'))\leq
Cd'(\sigma,\sigma')^\alpha.\]
\end{Lemma}

\begin{proof}  Applying the average contractiveness condition
$-m+1$ times gives
\begin{eqnarray*}
\sum\limits_{e_{m},...,e_0}&&\frac{1}{N}p_{e_{m}}(x_{i(e_{m})})...p_{e_{0}}(w_{e_{-1}}\circ...\circ
w_{e_{m}}(x_{i(e_{m})}))\\
&&\times d(w_{e_{0}}\circ...\circ
w_{e_{m}}(x_{i(e_{m})}),w_{e_{0}}\circ...\circ
w_{e_{m}}(y_{i(e_{m})}))\\
&\leq& a^{-m+1}\frac{1}{N}\sum\limits_{i=1}^Nd(x_i,y_i),
\end{eqnarray*}
 i.e.
\[\int d\left(Y^{x_1...x_N}_{m0},Y^{y_1...y_N}_{m0}\right)dP^m_{x_1...x_N}\leq a^{-m+1}
\frac{1}{N}\sum\limits_{i=1}^Nd(x_i,y_i).\] So, by Markov
inequality,
 \[P^m_{x_1...x_N}\left(d\left(Y^{x_1...x_N}_{mn},Y^{y_1...y_N}_{mn}\right)>a^{\frac{-m+1}{2}}
 \frac{1}{N}\sum\limits_{i=1}^Nd(x_i,y_i)\right)\leq a^{\frac{-m+1}{2}}.\]
Set $A_{m}:=\left\{\sigma\in\Sigma:\
d\left(Y^{x_1...x_N}_{m0}(\sigma),Y^{y_1...y_N}_{m0}(\sigma)\right)>a^{\frac{-m+1}{2}}
\frac{1}{N}\sum_{i=1}^Nd(x_i,y_i)\right\}$ and
 \[A:=\bigcap\limits_{l=0}^{-\infty}\bigcup\limits_{m=l}^{-\infty} A_{m}.\]
Then
\begin{eqnarray*}
 P_{x_1...x_N}(A)\leq P_{x_1...x_N}\left(\bigcup\limits_{m=l}^{-\infty}
 A_m\right)\leq \sum\limits_{m=l}^{-\infty}P^m_{x_1...x_N}(A_m)\leq
 \sum\limits_{m=l}^{-\infty}a^{\frac{-m+1}{2}},
\end{eqnarray*}
since $(\emptyset,...,\emptyset,
A_l,A_{l-1},...)\in\mathcal{C}\left(\bigcup_{m=l}^{-\infty}
 A_m\right)$ for all $l\leq 0$. Hence
$P_{x_1...x_N}(A)=0$ and for every $\sigma\in\Sigma\setminus A$
\[d\left(Y^{x_1...x_N}_{m0}(\sigma),Y^{y_1...y_N}_{m0}(\sigma)\right)\leq a^{\frac{-m+1}{2}}
\frac{1}{N}\sum\limits_{i=1}^Nd(x_i,y_i)\] for all $m$ except
finitely many. This implies $(i)$.

Now, for part $(ii)$ set $C:=\max_{e\in
E}d(x_{t(e)},w_e(x_{i(e)}))$. Then applying the average
contractiveness condition $-m+1$ times reveals that
\begin{eqnarray*}
\int
d\left(Y^{x_1...x_N}_{m0},Y^{x_1...x_N}_{(m-1)0}\right)dP^m_{x_1...x_N}&\leq&
a^{-m+1} \frac{1}{N}\sum\limits_{e\in
E}p_e(x_{i(e)})d(x_{t(e)},w_e(x_{i(e)}))\\
&\leq&a^{-m+1}C.
\end{eqnarray*}
So, by Markov inequality,
 \[P^m_{x_1...x_N}\left(d\left(Y^{x_1...x_N}_{m0},Y^{x_1...x_N}_{(m-1)0}\right)>a^{\frac{-m+1}{2}}
 C\right)\leq a^{\frac{-m+1}{2}}.\]
Set $B_{m}:=\left\{\sigma\in\Sigma:\
d\left(Y^{x_1...x_N}_{m0}(\sigma),Y^{x_1...x_N}_{(m-1)0}(\sigma)\right)>a^{\frac{-m+1}{2}}
C\right\}$ and
 \[B:=\bigcap\limits_{l=0}^{-\infty}\bigcup\limits_{m=l}^{-\infty} B_{m}.\]
Then
\begin{eqnarray*}
 P_{x_1...x_N}(B)\leq P_{x_1...x_N}\left(\bigcup\limits_{m=l}^{-\infty}
 B_m\right)\leq \sum\limits_{m=l}^{-\infty}P^m_{x_1...x_N}(B_m)\leq
 \sum\limits_{m=l}^{-\infty}a^{\frac{-m+1}{2}},
\end{eqnarray*}
since $(\emptyset,...,\emptyset,
B_l,B_{l-1},...)\in\mathcal{C}\left(\bigcup_{m=l}^{-\infty}
 B_m\right)$ for all $l\leq 0$.
Hence $P_{x_1...x_N}(B)=0$ and, for every
$\sigma\in\Sigma\setminus B$,
 \[\sum\limits_{m=0}^{-\infty}
d\left(Y^{x_1...x_N}_{m0}(\sigma),Y^{x_1...x_N}_{(m-1)0}(\sigma)\right)<\infty.\]
This implies that
$\left(Y^{x_1...x_N}_{m0}(\sigma)\right)_{m\in\mathbb{Z}\setminus\mathbb{N}}$
is a Cauchy sequence for $P_{x_1...x_N}$-a.e. $\sigma\in\Sigma$,
and so $F_{x_1...x_N}:=\lim_{m\to-\infty}Y^{x_1...x_N}_{m0}$
exists $P_{x_1...x_N}$-a.e..

(iii) Let
\[Q_l:=\Sigma\setminus\bigcup\limits_{m=l}^{-\infty}B_m\mbox{ for all } l\in\mathbb{Z}\setminus
\mathbb{N}.\]
 Then $Q_l\subset Q_{l-1}$
for all $l\in\mathbb{Z}\setminus\mathbb{N}$. Since every
$Y^{x_1...x_N}_{m0}$ is continuous, $B_m$ is open for all
$m\in\mathbb{Z}\setminus\mathbb{N}$, and so, every $Q_l$ is
closed. By the above,
\[\sum\limits_{l=0}^{-\infty} P_{x_1...x_N}(\Sigma\setminus Q_l)=\sum\limits_{l=0}^{-\infty}\sum\limits_{m=l}^
{-\infty}a^{\frac{-m+1}{2}}<\infty.\] Now, let $\sigma,\sigma'\in
Q_l$ for some $l$. Then
\[d\left(Y^{x_1...x_N}_{m0}(\sigma),Y^{x_1...x_N}_{(m-1)0}(\sigma)\right)\leq
a^{\frac{-m+1}{2}}C\mbox{ for all }m\leq l.\] Therefore, by the
triangle inequality,
\[d\left(Y^{x_1...x_N}_{m0}(\sigma),Y^{x_1...x_N}_{(m-k)0}(\sigma)\right)\leq\sum\limits_{i=m}^{-\infty}
a^{\frac{-i+1}{2}}C=C\frac{1}{1-\sqrt{a}}\ a^{\frac{|m-1|}{2}}\mbox{
for all }m\leq l,\ k\geq 1.\] Hence
\[d\left(Y^{x_1...x_N}_{m0}(\sigma),F_{x_1...x_N}(\sigma)\right)\leq
 C\frac{1}{1-\sqrt{a}}\
 a^{\frac{|m-1|}{2}}\mbox{ for all }m\leq l.\]
 Analogously,
 \[d\left(Y^{x_1...x_N}_{m0}(\sigma'),F_{x_1...x_N}(\sigma')\right)\leq
 C\frac{1}{1-\sqrt{a}}\
 a^{\frac{|m-1|}{2}}\mbox{ for all }m\leq l.\] Now, let $\delta_l=(1/2)^{|l-1|}$,
$\alpha:=\log\sqrt{a}/\log(1/2)$ and
$d'(\sigma,\sigma')=(1/2)^{|m-1|}$ for some $m\leq l$. Then
$Y^{x_1...x_N}_{m0}(\sigma')=Y^{x_1...x_N}_{m0}(\sigma)$.
Therefore,
\[d\left(F_{x_1...x_N}(\sigma),F_{x_1...x_N}(\sigma')\right)\leq
\frac{2C}{1-\sqrt{a}}\ a^\frac{|m-1|}{2}=\frac{2C}{1-\sqrt{a}}\
d'(\sigma,\sigma')^\alpha.\] \end{proof}

\begin{Definition}
We call the map
  \begin{eqnarray*}
    F_{x_1...x_N}:\Sigma\longrightarrow K,
 \end{eqnarray*}
which is defined  $P_{x_1...x_N}$-a.e. by Lemma \ref{cml},
 the {\it coding map} of the CMS.
\end{Definition}

\section{Definition with respect to a generalized Markov measure}

Our next aim is to show that the coding map is defined almost
everywhere with respect to any outer measure $\Phi(\nu)$ if each
probability function $p_e|_{K_{i(e)}}$ is Dini-continuous and
bounded away from zero by $\delta>0$. For that, we only need to
establish that $\Phi(\nu)$ is absolutely continuous with respect
to $P_{x_1...x_N}$ in this case.

\begin{Definition}
 We call a function $f:(X,d)\longrightarrow\mathbb{R}$  {\it Dini-continuous} iff
 for some $c >0$  \[\int_0^c\frac{\phi(t)}{t}dt<\infty\]
  where $\phi$ is {\it the modulus of uniform continuity} of $f$, i.e.
     \[\phi(t):=\sup\{|f(x)-f(y)|:d(x,y)\leq t,\ x,y\in X\}.\]
\end{Definition}
It is easily seen that the Dini-continuity is weaker than the
H\"{o}lder and stronger than the uniform continuity.  There  is a
well known characterization of the Dini-continuity, which will be
useful later.
\begin{Lemma}\label{Dc}
 Let $0<c<1$ and $b>0$.
  A  function $f$ is Dini-continuous  iff
 \[\sum_{n=0}^\infty\phi\left(bc^n\right)<\infty\] where $\phi$ is
 the modulus of uniform continuity of $f$.
\end{Lemma}
The proof is simple (e.g. see \cite{Wer1}).

Set $Z^x_{mn}(\sigma):=w_{\sigma_{n}}\circ...\circ
w_{\sigma_{m}}(x)$ for all $x\in K$, $m\leq n\in\mathbb{Z}$ and
$\sigma\in\Sigma$.
\begin{Lemma}\label{dd}
  Let $x_i \in K_i$ for every $1\leq i\leq N$ and $x\in K$. Let $i_0\in\{1,...,N\}$ such that $x\in K_{i_0}$.
  Then for all integers $m\leq 0$ and for all $\epsilon>0$ there exist $k\geq m$
and $B\in\mathcal{A}_m$ such that $P^m_{x}(B)<\epsilon$ and
\[ n\geq k\ \Rightarrow\ d\left(Z^x_{mn}(\sigma),Y^{x_1...x_N}_{mn}(\sigma)\right)\leq
a^{\frac{n-m+1}{2}} d(x,x_{i_0})\] for all
$\sigma\in\Sigma\setminus B$.
\end{Lemma}
\begin{proof} Fix $m\leq 0$. Applying the average contractiveness
condition $n-m+1$ times gives
\begin{eqnarray*}
&&\sum\limits_{e_{m},...,e_n}p_{e_{m}}(x)...p_{e_{n}}(w_{e_{n-1}}\circ...\circ
w_{e_{m}}(x))\\
&&\times  d(w_{e_{n}}\circ...\circ
w_{e_{m}}(x),w_{e_{n}}\circ...\circ w_{e_{m}}(x_{i(e_{m})})) \leq
a^{n-m+1}d\left(x,x_{i_0}\right),
\end{eqnarray*}
 i.e.
\[\int d\left(Z^x_{mn},Y^{x_1...x_N}_{mn}\right)dP^m_{x}\leq a^{n-m+1}
d(x,x_{i_0}).\] So, by the Markov inequality,
 \[P^m_{x}\left(d\left(Z^x_{mn},Y^{x_1...x_N}_{mn}\right)>a^{\frac{n-m+1}{2}}
 d(x,x_{i_0})\right)\leq a^{\frac{n-m+1}{2}}.\]
Set $A_{mn}:=\left\{\sigma\in\Sigma:\
d\left(Z^x_{mn}(\sigma),Y^{x_1...x_N}_{mn}(\sigma)\right)>a^{\frac{n-m+1}{2}}
d(x,x_{i_0})\right\}$ for $m\leq n$ and $B_k:=\bigcup_{n\geq
k}^\infty A_{mn}$ for $k\geq m$. Then
\[\sum\limits_{n=m}^\infty P^m_{x}(A_{mn})\leq\sum\limits_{l=1}^\infty a^{\frac{l}{2}}<\infty. \]
Therefore
\[P^m_{x}\left(\bigcap\limits_{k=m}^\infty B_k\right)=0,\]
by the Borel-Cantelli argument. Hence, for all $\epsilon>0$ there
exists $k\geq m$ such that $P^m_{x}(B_k)<\epsilon$ and
\[n\geq k\ \Rightarrow\ d\left(Z^x_{mn}(\sigma),Y^{x_1...x_N}_{mn}(\sigma)\right)\leq
a^{\frac{n-m+1}{2}} d(x,x_{i_0})\] for all
$\sigma\in\Sigma\setminus B_k$.\end{proof}

The next lemma is a generalization of Lemma 3 in \cite{Elton}. The
proof of it that of Elton.
\begin{Lemma}\label{abc}
Suppose that  $p_e|_{K_{i(e)}}$ is Dini-continuous and there
exists $\delta>0$ such that $p_e|_{K_{i(e)}}\geq\delta$ for all
$e\in E$.
 Let $x_i\in K_i$ for all  $1\leq i\leq N$ and $x\in K$.
  Then $P^m_x$ is absolutely continuous with respect to $P^m_{x_1...x_N}$  for all $m\leq 0$.
\end{Lemma}

\begin{proof} Fix $m\leq 0$. Let $A\in\mathcal{A}_m$ be such that
$P^m_{x_1...x_N}(A)=0$ and $\epsilon>0$. We show
$P^m_x(A)<\epsilon$.

Let $i_0\in\{1,...,N\}$ such that $x\in K_{i_0}$. By Lemma
\ref{dd}, there exists $n_\epsilon\geq m$ and $B\in\mathcal{A}_m$
such that $P^m_x(B)<\epsilon/2$ and
\[ n\geq n_\epsilon\ \Rightarrow\ d\left(Z^x_{mn}(\sigma),Y^{x_1...x_N}
_{mn}(\sigma)\right)\leq a^{\frac{n-m+1}{2}}d(x,x_{i_0})\] for all
$\sigma\in\Sigma\setminus B$. Let $\phi_e$ be the modulus of
uniform continuity of $p_e|_{K_{i(e)}}$ for each $e\in E$ and
$\phi:=\max_{e\in E}\phi_e$. Since each $p_e|_{K_{i(e)}}$ is
Dini-continuous, by Lemma \ref{Dc}, we can choose $l\geq
n_\epsilon$ such that
$\sum_{k=l+1}^\infty\phi\left(a^{(k-m)/2}d(x,x_{i_0})\right)<\delta/2$.
Set
\begin{eqnarray*}Q_{n}:= \left\{
\begin{array}{l}\left\{\sigma\in\Sigma:  d(Z^x_{mk}(\sigma),Y^{x_1...x_N}_{mk}(\sigma)) \leq
a^{\frac{k-m+1}{2}}d(x,x_{i_0})\ \forall l\leq k\leq n\right\},\\ \mbox{ if }l\leq n\\
\Sigma, \mbox{ else }\end{array}\right.
\end{eqnarray*}
and $Q:=\bigcap_{n\geq m}Q_{n}$. Then $\Sigma\setminus B\subset Q$
and therefore $P_x(\Sigma\setminus Q)<\epsilon/2$. Now, for
$\sigma\in Q_{n}$, if $l\leq n$ and $(\sigma_m,...,\sigma_n)$ is a
path of the digraph starting in $i_0$, then
\begin{eqnarray*}
& &p_{\sigma_m}(x)...p_{\sigma_n}(w_{\sigma_{n-1}}\circ...\circ
w_{\sigma_m}x)\\
&\leq&p_{\sigma_m}(x_{i(\sigma_m)})...p_{\sigma_n}(w_{\sigma_{n-1}}\circ...\circ
w_{\sigma_m}x_{i(\sigma_m)})\left(\frac{1-\delta}{\delta}\right)^{l-m+1}\\
&\times&\prod\limits_{k=l+1}^n\left[1+\frac{p_{\sigma_k}(w_{\sigma_{k-1}}\circ...\circ
w_{\sigma_{m}}x)-p_{\sigma_k}(w_{\sigma_{k-1}}\circ...\circ
w_{\sigma_{m}}x_{i(\sigma_m)})}{p_{\sigma_k}(w_{\sigma_{k-1}}\circ...\circ w_{\sigma_{m}}x_{i(\sigma_m)})}\right]\\
&\leq&p_{\sigma_m}(x_{i(\sigma_m)})...p_{\sigma_n}(w_{\sigma_{n-1}}\circ...\circ
w_{\sigma_m}x_{i(\sigma_m)})\left(\frac{1-\delta}{\delta}\right)^{l-m+1}\\
&\times&\prod\limits_{k=l+1}^n\left[1+\frac{\phi\left(a^{\frac{k-m}{2}}d(x,x_{i_0})\right)}{\delta}\right].
\end{eqnarray*}
Since
$\prod_{k=l+1}^n\left[1+{\phi\left(a^{(k-m)/{2}}d(x,x_{i_0})\right)}/{\delta}\right]\leq
1+2\sum_{k=l+1}^\infty{\phi\left(a^{(k-m)/2}d(x,x_{i_0})\right)}/{\delta}\leq2$,
it follows that
\begin{eqnarray*}
 p_{\sigma_m}(x)...p_{\sigma_n}(w_{\sigma_{n-1}}\circ...\circ
w_{\sigma_m}x)&\leq&
2N\left(\frac{1-\delta}{\delta}\right)^{l-m+1}\frac{1}{N}\\
&&\times
p_{\sigma_m}(x_{i(\sigma_m)})...p_{\sigma_n}(w_{\sigma_{n-1}}\circ...\circ
w_{\sigma_m}x_{i(\sigma_m)}).
\end{eqnarray*}
 If  $l> n$ or
$(\sigma_m,...,\sigma_n)$ is not a path of the digraph starting in
$i_0$, then it holds trivially for any $\sigma\in\Sigma$.

Let $\triangle_m$ be the algebra every element of which is a finite
union of cylinders of the form $_m[e_m,...,e_n]$, $n\geq m$. By
Caratheodory construction, there exists a sequence
$(A_k)_{k\in\mathbb{N}}\subset\triangle_m$ such that
$A\subset\bigcup_{k=1}^\infty A_k$ and
\[\sum\limits_{k=1}^\infty
P^m_{x_1...x_N}(A_k)<\frac{\epsilon}{4N}\left(\frac{\delta}{1-\delta}\right)^{l-m+1}.\]
 We can write every finite union $\bigcup_{k=1}^nA_k$ as a disjoint
union $\biguplus_{k=1}^{m_n} C_k$ of cylinders which generate
$\triangle_m$. Let $ C_k=_m[e_m,...,e_n]$ with $m\leq n$. Then
\begin{eqnarray*}
  P^m_x\left(Q\cap C_k\right)&\leq& P^m_x\left(Q_{n}\cap
  C_k\right)\\
  &=&\sum\limits_{(\sigma_m,...,\sigma_n):\sigma\in Q_{n}\cap
  C_k}p_{\sigma_m}(x)...p_{\sigma_n}(w_{\sigma_{n-1}}\circ...\circ
  w_{\sigma_m}x)\\
  &\leq& 2N\left(\frac{1-\delta}{\delta}\right)^{l-m+1}\sum\limits_{(\sigma_m,...,\sigma_n): \sigma\in Q_{n}\cap
  C_k}\frac{1}{N}\\
  &&\times p_{\sigma_m}(x_{i(\sigma_m)})...p_{\sigma_n}(w_{\sigma_{n-1}}\circ...\circ
w_{\sigma_m}x_{i(\sigma_m)})\\
  &\leq& 2N\left(\frac{1-\delta}{\delta}\right)^{l-m+1} P^m_{x_1...x_n}\left( C_k\right).
\end{eqnarray*}
Hence
\begin{eqnarray*}
  P^m_x(A)&=&P^m_x(Q\cap A)+P^m_x(A\setminus Q)\\
  &\leq&\lim\limits_{n\to\infty}P^m_x\left(\biguplus_{k=1}^{m_n}
  C_k\cap Q\right)+\frac{\epsilon}{2}\\
  &=&\lim\limits_{n\to\infty}\sum\limits_{k=1}^{m_n}P^m_x\left(
  C_k\cap Q\right)+\frac{\epsilon}{2}\\
  &\leq& 2N\lim\limits_{n\to\infty}\left(\frac{1-\delta}{\delta}\right)^{l-m+1}\sum\limits_{k=1}^{m_n} P^m_{x_1...x_n}\left(
  C_k\right)+\frac{\epsilon}{2}\\
  &\leq& 2N\left(\frac{1-\delta}{\delta}\right)^{l-m+1}\sum\limits_{k=1}^{\infty} P^m_{x_1...x_n}\left(
  A_k\right)+\frac{\epsilon}{2}\\
  &<&\epsilon.
\end{eqnarray*}
\end{proof}
\begin{Theorem}\label{cmT}
Suppose that  $p_e|_{K_{i(e)}}$ is Dini-continuous and there
exists $\delta>0$ such that $p_e|_{K_{i(e)}}\geq\delta$ for all
$e\in E$.
 Let $x_i,y_i\in K_i$ for all  $1\leq i\leq N$ and $\nu\in P(K)$. Then:

 (i)  $F_{x_1...x_N}$ is defined $\Phi(\nu)$-a.e.,

 (ii) $F_{x_1...x_N}=F_{y_1...y_N}$ $\Phi(\nu)$-a.e.. and

 (iii) There exists a sequence of closed subsets $Q_1\subset
 Q_2\subset...\subset\Sigma$ with\\ $\lim_{k\to\infty}\Phi(\nu)(\Sigma\setminus Q_k)=0$
 such that all $F_{x_1...x_N}|_{Q_k}$ are locally
 H\"{o}lder-continuous with the same H\"{o}lder-constants.
\end{Theorem}

\begin{proof} By Lemma \ref{abc}, $\Phi_m(\nu)$ is absolutely
continuous with respect to $P^m_{x_1...x_N}$ for all $m\leq 0$. By
Lemma \ref{om} $(ii)$, this implies that $\Phi(\nu)$ is absolutely
continuous with respect to $P_{x_1...x_N}$. The claim follows by
Lemma \ref{cml}.\end{proof}

\begin{Definition}
We call
\[M:=\Phi(\mu)\]
a {\it generalized Markov measure}, where $\mu$ is the invariant
Borel probability measure of the CMS. Denote the Borel
$\sigma$-algebra on $\Sigma$ by $\mathcal{B}(\Sigma)$.
\end{Definition}

\begin{Proposition}
 $M$ is a shift invariant Borel probability
measure on $\Sigma$ with
\[M\left(_m[e_1,...,e_k]\right):=\int
p_{e_1}(x)p_{e_2}(w_{e_1}x)...p_{e_k}(w_{e_{k-1}}\circ...\circ
w_{e_1}x)d\mu(x)\] for every cylinder set
$_m[e_1,...,e_k]\subset\Sigma.$
\end{Proposition}

\begin{proof} First, define a set function $\bar M$ on all
cylinders of $\Sigma$ by
\[\bar M( _m[e_m,...,e_n]):=\Phi_m(\mu)(_m[e_m,...,e_n]).\]
 We show that $\bar M$ extends uniquely to a shift invariant Borel probability measure on $\Sigma$ and
 \[\bar M|_{\mathcal{A}_m}=\Phi_m(\mu)\mbox{ for all }m\leq 0.\]
 We only need to check that
\begin{eqnarray*}
\Bar M( _m[e_m,...,e_n])=\sum\limits_{e_{n+1}}\Bar M(
_m[e_m,...,e_n,e_{n+1}])
\end{eqnarray*}
and that
\begin{eqnarray*}
\Bar M( _m[e_m,...,e_n])=\sum\limits_{e_{m-1}}\bar M(
_{m-1}[e_{m-1},e_m,...,e_n]),
\end{eqnarray*}
the rest follows by the standard extension argument. The first
equation is obvious by the definition of $P_x$. For the second we
need the invariance of $\mu$.
\begin{eqnarray*}
& &\sum\limits_{e_{m-1}}\bar M( _{m-1}[e_{m-1},e_m,...,e_n])\\
&=&\sum\limits_{e_{m-1}}\int
p_{e_{m-1}}(x)p_{e_m}(w_{e_{m-1}}(x))...p_{e_n}(w_{e_{n-1}}\circ...\circ
w_{e_{m-1}}(x))d\mu(x)\\
&=&U^*\mu(p_{e_m}...p_{e_n}\circ
w_{e_{n-1}}\circ...\circ w_{e_{m}})\\
&=&\mu(p_{e_m}...p_{e_n}\circ
w_{e_{n-1}}\circ...\circ w_{e_{m}})\\
&=&\bar M( _m[e_m,...,e_n]).
\end{eqnarray*}
 Now, we show that
 \[\bar M=M|_{\mathcal{B}(\Sigma)}.\]
 Let $B\in \mathcal{B}(\Sigma)$ and
$\epsilon>0$. Since $\mathcal{B}(\Sigma)$ is the smallest
$\sigma$-algebra containing all $\mathcal{A}_m$, $m\leq 0$, it
follows that for every $(A_m)_{m\leq 0}\in\mathcal{C}(B)$
\begin{eqnarray*}
\bar M(B)\leq \bar M\left(\bigcup\limits_{m\leq
  0}A_m\right)\leq\sum\limits_{m\leq 0}\bar M(A_m)=\sum\limits_{m\leq 0} \Phi_m(\mu)(A_m).
\end{eqnarray*}
Hence $\bar M(B)\leq M (B)$. On the other hand, let $\triangle$ be
the algebra generated by all cylinders in $\Sigma$. Then every
$A\in\triangle$ is also an element of some $\mathcal{A}_m$, $m\leq
0$, and therefore $\bar M(A)\geq M (A)$.  By Caratheodory
construction, there exists  a sequence
$(A_k)_{k\in\mathbb{N}}\subset\triangle$ such that
$B\subset\bigcup_{k\in\mathbb{N}}A_k$ and
\[\bar M(B)>\sum\limits_{k\in\mathbb{N}}\bar M(A_k)-\epsilon.\]
Therefore
\[\bar M(B)\geq\sum\limits_{k\in\mathbb{N}}
M(A_k)-\epsilon\geq
M\left(\bigcup\limits_{k\in\mathbb{N}}A_k\right)-\epsilon\geq
M(B)-\epsilon.\] Hence, $\bar M(B)= M (B)$.\end{proof}

\begin{Definition}\label{gMs}
We call the measure preserving transformation $S$ of the
probability space $\left(\Sigma, \mathcal{B}(\Sigma), M\right)$ a
{\it generalized Markov shift}.
\end{Definition}

  Now, we state explicitly the important special case of Theorem \ref{cmT}.
\begin{Corollary}\label{cmC}
Suppose that  $p_e|_{K_{i(e)}}$ is Dini-continuous and there
exists $\delta>0$ such that $p_e|_{K_{i(e)}}\geq\delta$ for all
$e\in E$.
 Let $x_i,y_i\in K_i$ for all  $1\leq i\leq N$. Then:

 (i)  $F_{x_1...x_N}$ is defined $M$-a.e.,

 (ii) $F_{x_1...x_N}=F_{y_1...y_N}$ $M$-a.e.. and

 (iii) There exists a sequence of closed subsets $Q_1\subset
 Q_2\subset...\subset\Sigma$ with\\ $\lim_{k\to\infty}M(Q_k)=1$
 such that all $F_{x_1...x_N}|_{Q_k}$ are locally
 H\"{o}lder-continuous with the same H\"{o}lder-constants.
\end{Corollary}

\begin{acknowledgements}
I would like to thank: EPSRC and School of Mathematics
 and Statistics of University of St Andrews for providing me with a scholarship and excellent working
 conditions in St Andrews, Professor K. J. Falconer for valuable comments on the first draft of the paper,
 my supervisor Lars Olsen for his interest
 in my work, his support and many fruitful discussions, and the
 anonymous referee for various improvements to this paper.
\end{acknowledgements}

\end{document}